\RequirePackage{etex}

\documentclass[AMS,STIX1COL]{WileyNJD-v2}

\articletype{\emph{author's accepted manuscript by Wiley publisher}}

\begin{document}

\color{black}

\title{Nonlinear integral extension of PID control with improved
convergence of perturbed second-order dynamic systems}

\author{Michael Ruderman*}

%\author[2,3]{Author Two}

%\author[3]{Author Three}

\authormark{MICHAEL RUDERMAN}

\address{\orgdiv{Department of Engineering Sciences}, \orgname{University of Agder}, \orgaddress{\state{P.B. 422, Kristiansand, 4604}, \country{Norway}}}

%\address[2]{\orgdiv{Org Division}, \orgname{Org Name}, \orgaddress{\state{State name}, \country{Country name}}}

%\address[3]{\orgdiv{Org Division}, \orgname{Org Name}, \orgaddress{\state{State name}, \country{Country name}}}

\corres{* \email{michael.ruderman@uia.no}}

%\presentaddress{This is sample for present address text this is
%sample for present address text}

%%%
\abstract[Summary]{Nonlinear extension of the integral part of a
standard proportional-integral-derivative (PID) feedback control
is proposed for perturbed second-order systems. The approach is model-free and requires solely the Lipschitz boundedness of the unknown matched perturbations. For
constant disturbances, the global asymptotic stability is
shown based on the circle criterion. For Lipschitz perturbations, an ultimately bounded
output error is provided based on the steady-state behavior in frequency domain.
Also the transient response to the stepwise disturbances is analyzed for the control tuning. Based on the developed analysis, the design recommendations are formulated as a step by step procedure. It is also discussed how the proposed control
is applicable to second-order systems extended by additional
(parasitic) actuator dynamics with low-pass characteristics. The
proposed nonlinear control is proven to outperform its linear PID
counterpart during the settling phase, i.e. at convergence of the
residual output error. An experimental case study of the
second-order system with an additional actuator dynamics and
considerable perturbations is demonstrated to confirm and
benchmark the control performance.}


\maketitle

%\footnotetext{\textbf{Abbreviations:} ANA, anti-nuclear
%antibodies; APC, antigen-presenting cells; IRF, interferon
%regulatory factor}

\newtheorem{thm}{Theorem}

%\newtheorem{lem}[thm]{Lemma}
%\newtheorem{clr}{Corollary}
%\newdefinition{remark}{Remark}
%\newdefinition{exmp}{Example}
%\newdefinition{prop}{Proposition}
%\newproof{pf}{Proof}

%\newtheorem{theorem}{Theorem}
\newtheorem{rem}{Remark}

%%%%%%%%%%%%%%%%%%%%%%%%%%%%%%%%%%%%%%%%%%%%%%%%%%%%%%%%%%%%%%%%%%%%%%%%%%%%%%%%
\section{INTRODUCTION}
\label{sec:1}

Proportional-integral-derivative (PID) controllers, see e.g.
\cite{aastrom2006,aastrom2021}, as well as their possible 
extensions like anti-windup \cite{hippe2006}, gain-scheduling,
parameters adaptation, \cite{gyongy2006automatic}, and others
\cite{ortega2021pid} are still seen as a mostly working horse in
various industries. Despite a large number of the available and
well perceived tutorials and surveys on the PID principles, tuning,
extensions, and limitations, e.g.
\cite{ang2005pid,hara2006robust,garpinger2014performance} just for
mentioning here the few, numerous applications face here and again
some crucial challenges with use of the PID controllers. This is
often due to a trade-off between the respondness and robustness, while
the performance metrics associated with a transient and
steady-state response appear to be most significant for several
applications.

Well known, an integral control action is required for asymptotic
set-point regulation under parameter perturbations or in presence of
the matched stationary disturbances, which are quite common in
the control practice. For simple first- and second-order
processes, including also the possible time delays, the tuning of
the integral time constant is straightforward and well understood,
see e.g. \cite{skogestad2003}. More challenging becomes an
effective integral control action in case of disturbance processes
that vary over time, especially in vicinity to the set-point, and
when the application requires an increased settling accuracy in
a shorter time. Hence it is not surprising that a standard integral
control action cannot compensate accurately, for example, the
nonlinear Coulomb friction \cite{bisoffi2017,ruderman2022stick,ruderman2025} in
the set-point tasks. As one of the possible solutions, the
so-called reset integral control, which was introduced originally
in \cite{clegg1958}, can be applied, see e.g. \cite{beerens2019}.
The reset integral behavior in feedback, and thus a switching
dynamics, restricts however the applicability to a certain class
of the system plants \cite{beker2004}. Worth noting is also that
despite a relaxation applied to the underlying base systems and,
this way, an attempt to address the stability by a less stringent
circle-criterion conditions \cite{van2017}, the reset integrator
does not entirely satisfy the necessary conditions posed on the
feedback nonlinearity, cf. \cite{khalil2002,sastry2013}, cf. with
preliminaries provided below in section \ref{sec:2}. Another
approach of the time-varying integral control gain, which assumes
a monotonic decrease of the gain to zero at a sufficiently slow
rate, was presented in \cite{logemann2000} for linear system
plants subject to a globally Lipschitz and non-decreasing actuator
nonlinearity. Despite the profound theoretical results, which also
form the basis for a simple adaptive control strategy, the approach is
not directly applicable to compensate for external disturbances
and to improve the convergence performance. An universal integral
controller for minimum-phase nonlinear systems and the tuning
procedures for its key parameters were proposed in
\cite{khalil2000}. A single-input-single-output (SISO) nonlinear
system plant that has a well-defined normal form with
asymptotically stable zero dynamics was assumed for the synthesis
procedure. Remarkable is that for relative degree-one systems,
this universal controller reduces to the classical PI controller
followed by the saturation, and for relative degree-two systems it
reduces to the classical PID controller followed by the saturation.
Settling performance, in terms of the convergence time and
accuracy were, however, not explicitly in focus of the proposed
methodology. Later, the problem of robust nonlinear integral
control for a class of systems with dynamic uncertainties was
addressed in \cite{jiang2001}, while a constructive controller
design procedure was presented in the presence of an unknown
equilibrium due to uncertain nonlinearities and dynamic
uncertainties. The resulted controller reduces, notwithstanding,
to linear control laws if the nominal system is linear, and
becomes a traditional PI control if, additionally, the relative
degree is one, or the uncertainty falls into the input space, i.e.
it fulfills the strict matching condition.

Motivated by the background provided above and the perturbed steady-state
behavior of motion control systems, see e.g.
\cite{ruderman2022motion}, the present work introduces a possible
nonlinear integral extension of the PID feedback control. To the
best of the author's knowledge, the proposed control, see eq.
\eqref{eq:3:1}, is new and was not reported in the previously
published works. An innovative contribution of the work lies in extending a standard linear PID feedback controller by a nonlinear integral part which allows for faster output convergence in presence of the unknown matched Lipschitz perturbations. Such extension is model-free and relies solely on the Lipschitz continuity of the disturbing perturbations. A further contribution is in the simple analysis of stability and output error dynamics as well as in the illustrative experimental case study that supports the proposed control approach.

The rest of the paper is structured in the following way. Section
\ref{sec:2} provides the notation in use and the necessary
preliminaries of a standard PID feedback control and the circle
criterion of the absolute stability. The main results are placed into
section \ref{sec:3}, starting with absolute stability, followed by
the disturbance rejection and extension of the second-order system
plant by an additional actuator dynamics. Also some clarifying
numerical examples are included, and a step-by-step controller design procedure is recommended. An experimental case study used for evaluation of the proposed control and benchmarking it against the standard PD and PID controllers is provided in section
\ref{sec:4}. The paper ends with the short conclusions drawn in section \ref{sec:5}.

%%%%%%%%%%%%%%%%%%%%%%%%%%%%%%%%%%%%%%%%%%%%%%%%%%%%%%%%%%%%%%%%%%%%%%%%%%%%%%%%
\section{PRELIMINARIES}
\label{sec:2}

\subsection*{Notations}
\label{sec:2:sub:0}

Unless otherwise said, $x = [x_1, x_2, \ldots, x_n]^\top \in
\mathbb{R}^n $ is a state vector of real numbers, $t \in
\mathbb{R}_{+} = \{ \tau \in \mathbb{R} \: : \: \tau \geq 0 \}$ is
a non-negative time variable. The time derivative of a variable
$x_n$ is denoted by $\dot{x}_n \equiv d x_n / dt$, respectively
the second time derivative is denoted by $\ddot{x}_n \equiv d^2
x_n / dt^2 $, and so on. An input-output transfer function is
denoted by $H(s)=y(s)/u(s)$, where $s = (\rho + j \omega) \in
\mathbb{C}$ is the complex Laplace variable with the signed
angular frequency $\omega \in \mathbb{R}$, imaginary unit $j^2 =
-1$, and $\rho \in \mathbb{R}$. An absolute value in $\mathbb{R}$
is denoted by $|\cdot|$, while $\bigl|H(s)\bigr|$ and respectively
$\bigl|H(j\omega)\bigr|$ denote the absolute value of a complex
function. An identity matrix of the appropriate dimension is
denoted by $I$.

\subsection{PID controlled feedback loop}
\label{sec:2:sub:1}

We consider the second-order systems of the type
\begin{equation}\label{eq:2:1}
\ddot{y}(t) = u(t) + \sigma(t),
\end{equation}
where the output of interest $y(t)$ and its time derivative
$\dot{y}(t)$ are available for a feedback control. The control
signal $u(t) \in \mathbb{R}$ is continuous on $t \in
\mathbb{R}_{+}$ and assumed to stay unsaturated\footnote{We notice
that a feedback control involving an integral action may be
subject to saturation, that is quite common in practice of
the control engineering. However, an explicit measure against the
associated wind-up effects is beyond the scope of the current
work. Therefore, neither standard PID control, assumed as a
reference benchmarking control system, nor the proposed nonlinear
extension handle a saturated operation mode.} for $(y,\dot{y}) \in
\mathbb{R}^2$. The matched perturbation function $\sigma(t)$ is
assumed to be Lipschitz, i.e. $|\dot{\sigma}| < \Sigma$ for some
positive known $\Sigma$.

First, we recall briefly a standard PID feedback control, cf. e.g. \cite{aastrom2006,aastrom2021}, which in case of zero reference has the form 
\begin{equation}\label{eq:2:2}
u(t) = - a \, \dot{y}(t) - b \, y(t) - c \int y(t) dt,
\end{equation}
with the design parameters $a, b, c > 0$. Below we are focusing on
the most relevant (for us) properties and features of the PID control, while the control
parameters are assumed to be appropriately determined.

One can easily show, either in time or frequency domain, that for
a constant disturbance, i.e. $\dot{\sigma} = 0$, the PID feedback
control \eqref{eq:2:2} provides a global asymptotic stabilization
of \eqref{eq:2:1}. For a Lipschitz disturbance, mapped into
frequency domain by $\Sigma \propto \omega$, the PID feedback
control \eqref{eq:2:2} provides an ultimately bounded output error
$|y(t)| < \varepsilon(\omega)$ for $t \rightarrow \infty$.
Regarding this statement, which is rather self-understood in the
linear feedback control theory, cf. e.g.
\cite{franklin2019,aastrom2021}, we will have a closer look below
in section \ref{sec:3:sub:2}.

Substituting \eqref{eq:2:2} into \eqref{eq:2:1} and analyzing the
characteristic polynomial of the unperturbed system dynamics, i.e.
\begin{equation}\label{eq:2:3}
s^3 + a s^2 + b s + c = 0,
\end{equation}
one can show by, for example Routh-Hurwitz stability criterion,
e.g. \cite{franklin2019}, that the closed-loop control system is
asymptotically (and exponentially) stable if and only if
\begin{equation}\label{eq:2:4}
a b - c > 0.
\end{equation}
That means every solution of the state trajectories is uniformly
bounded and satisfies $y(t) \rightarrow 0 $, $\dot{y}(t)
\rightarrow 0$ as $t \rightarrow \infty$.

\subsection{Circle criterion}
\label{sec:2:sub:2}

For a linear time invariant (LTI) system in the minimal
state-space realization, see e.g. \cite{antsaklis2007} for detail,
a nonlinearity $\phi(\cdot, \cdot)$ that is acting in feedback
and satisfying
\begin{eqnarray*}
  \phi(0, t)   &=& 0 \quad \forall \; t \in \mathbb{R}_{+}, \\
  z\, \phi(z, t) & \geq & 0 \quad \forall \; z \in \mathbb{R}, \: t \in
\mathbb{R}_{+}
\end{eqnarray*}
results in the so called Lur'e problem \cite{lur1944}. The latter
is also refereed to as the \emph{absolute stability problem} and
is formulated for a generic class of the systems given by
\begin{eqnarray}
\label{eq:2:5}
\nonumber  \dot{x}(t) &=& A x(t) + B u(t), \\
            z(t) &=& C x(t) + D u(t), \\
\nonumber   u(t)  &=& - \phi\bigl(z(t), t\bigr).
\end{eqnarray}
Note that the SISO systems are considered, so that $u, \, z \in
\mathbb{R}$, while $A, \, B, \,C, \,D$ are the matrices
(correspondingly vectors) of an appropriate dimension, and they contain
the system parameters. Here and for the rest of the section, we
will closely follow the notations and developments of the
\emph{circle criterion} \cite{narendra1964,sandberg1964} as it is
provided in \cite{sastry2013}; for additional explanations see
also e.g. \cite{slotine1991}.

For applying the circle criterion of absolute stability to
\eqref{eq:2:5}, the following necessary assumptions are placed.
\begin{enumerate}[i)]
    \item $A$, $B$, and $C$ are the minimal realization of a
    linear subsystem in \eqref{eq:2:5}, i.e. the pair $(A,B)$ is
    fully controllable and the pair $(A,C)$ is fully observable.
    \item The system matrix $A$ has all eigenvalues $\lambda_i$ lying in
    $\mathbb{C}_{-}^{\circ}$, i.e. $\mathrm{Re} \bigl[ \lambda_i (A) \bigr] < 0$ for $i = 1, \ldots,
    n$.
    \item The nonlinearity $\phi(\cdot,t)$ belongs to the sector
    $[k_1,k_2]$ for $k_2 \geq k_1 \geq 0$. That means
    $$
    k_1 z \leq \phi(z,t) \leq k_2 z.
    $$
\end{enumerate}
Moreover, it is assumed that
\begin{enumerate}[i)]
    \item[iv)] $\phi(z,t)$ is memoryless, possibly time-varying,
and piecewise continuous in $t$ and locally Lipschitz in $z$.
\end{enumerate}
For the input-output transfer function
\begin{equation}\label{eq:2:6}
H(s) = C (sI - A)^{-1} + D
\end{equation}
of the linear subsystem in \eqref{eq:2:5}, the absolute stability
of the closed-loop \eqref{eq:2:5} with the nonlinearity $\phi(\cdot,t)$ belonging
to the sector $[0,k]$, is guaranteed by the circle criterion
stated in the following theorem.

\vspace{1mm}

\begin{thm}
\label{theorem:1} Given a feedback system described by the
equations \eqref{eq:2:5}. Let the assumptions (i)-(iv) be
effective. Then, the origin of \eqref{eq:2:5} is globally
asymptotically stable if
\begin{equation}\label{eq:2:7}
\mathrm{Re} \bigl[ 1 + k  H(j\omega)\bigr] > 0, \quad \forall \:
\omega.
\end{equation}
\end{thm}
\vspace{2mm} The proof of Theorem \ref{theorem:1} can be found in
\cite[page~241]{sastry2013}.

\vspace{1mm}

It is worth recalling that the circle criterion given by
\eqref{eq:2:7} for the sector $[0,k]$ is equivalent to,
correspondingly can be transformed into (and also from), the
circle criterion formulated for the sector $[k_1,k_2]$ in terms of
a disc $\mathcal{D}(k_1, k_2)$, which is drawn in the polar
coordinates of the $H(j\omega)$-locus, see e.g.
\cite{slotine1991,sastry2013}.

%%%%%%%%%%%%%%%%%%%%%%%%%%%%%%%%%%%%%%%%%%%%%%%%%%%%%%%%%%%%%%%%%%%%%%%%%%%%%%%%
\section{MAIN RESULTS}
\label{sec:3}

The proposed feedback control law with a nonlinear extension of
the integral part is given by
\begin{equation}\label{eq:3:1}
u(t) = - a \, \dot{y}(t) - b \, y(t) - c \, \Bigl( 1 + d \exp
\bigl(e |y(t)|\bigr) \Bigr) \int y(t) dt,
\end{equation}
cf. \eqref{eq:2:2}. The linear control parameters $a,b,c > 0$ are in the same sense as
for \eqref{eq:2:2}. The nonlinear parameters $d \geq 0$ and $e <
0$ provide a smooth variation of the overall integral gain
\begin{equation}\label{eq:3:2}
\Omega \equiv c \, \Bigl( 1 + d \exp \bigl(e |y(t)|\bigr) \Bigr),
\end{equation}
that is monotonically increasing with decrease of the magnitude of
the control error. It is worth emphasizing that:
\begin{enumerate}[(i)]
    \item $c \leq \Omega \leq c+cd$;
    \item $\Omega \rightarrow c$ for a continuously growing $|y|$;
    \item the control \eqref{eq:3:1} recovers completely to the
benchmarking reference controller \eqref{eq:2:2} if $d=0$.
\end{enumerate}

In the following, we are first proving the absolute stability of an
unperturbed system \eqref{eq:2:1}, \eqref{eq:3:1} by using the
circle criterion, cf. with section \ref{sec:2:sub:2}. Afterwards,
the closed-loop system \eqref{eq:2:1}, \eqref{eq:3:1} with an
exogenous input $\sigma$ is analyzed in terms of the convergence
performance of $y(t)$, and that in comparison to the control law
\eqref{eq:2:2}, i.e. for $d = 0$ in \eqref{eq:3:1}. An illustrative numerical example discloses the output convergence properties in presence of a constant disturbance when varying $d$-amplitude. Here we recall that the integral control part (of a standard PID control) is expected to compensate for constant unknown disturbances at steady-state. Then, a system extension by the additional (parasitic) actuator dynamics is also addressed. Finally, we draw several recommendation lines for design (in terms of the parameters tuning) of the proposed nonlinear extended PID controller.

\subsection{Absolute stability}
\label{sec:3:sub:1}

The unperturbed closed-loop system \eqref{eq:2:1}, \eqref{eq:3:1}
is transferable into the Lur'e problem \eqref{eq:2:5} with $D=0$ as
\begin{equation}\label{eq:3:1:1}
A = \left(%
\begin{array}{ccc}
  0 & 1 & 0 \\
  0 & 0 & 1 \\
  -c & -b & -a \\
\end{array}%
\right), \,
B = \left(%
\begin{array}{c}
  0 \\
  0 \\
  1 \\
\end{array}%
\right) , \, C^\top =
\left(%
\begin{array}{c}
  1 \\
  0 \\
  0 \\
\end{array}%
\right),
\end{equation}
and $x = \bigl[\int y dt, y, \dot{y}\bigr]^\top$. The necessary
assumptions i) and ii), cf. with section \ref{sec:2:sub:2}, are
fulfilled since the state-space realization \eqref{eq:3:1:1} is
both fully controllable and observable, and one can guarantee that
the system matrix $A$ is Hurwitz under the condition that the
inequality \eqref{eq:2:4} is satisfied. The feedback nonlinearity
\begin{equation}\label{eq:3:1:2}
\phi(z,t) = c \, d \, \exp \bigl(e |x_2| \bigr) \, z
\end{equation}
proves to satisfy the necessary assumption iii), since it belongs
to the sector $[0, k]$ with $k = cd > 0$ for all $|x_2(t)| \in
\mathbb{R}_{+}$. Fulfillment of the necessary assumption iv)
follows from the fact that \eqref{eq:3:1:2} is single-valued with
respect to a $(z, |x_2|)$ pair, i.e. it is memoryless, and
continuous in $t$ and Lipschitz in $z$. Two of the last statements
are true due to the exponential term, and since $x_2(t) \in
\mathcal{C}^3$ for all $u(t)$ and $\sigma(t)$ to be Lipschitz.
Now, we are in the position to show the absolute stability of the
unperturbed control system \eqref{eq:2:1}, \eqref{eq:3:1}.

\vspace{1mm}

\begin{thm}
\label{theorem:2} Given an unperturbed feedback system
\eqref{eq:2:1}, \eqref{eq:3:1}. Then, its origin is globally
asymptotically stable if
\begin{equation}\label{eq:3:1:3}
\frac{c^2 d - a c d \, \omega^2}{(b \omega - \omega^3)^2 + (c - a
\, \omega^2)^2 } > -1, \quad \forall \: \omega.
\end{equation}
\end{thm}

\vspace{1mm}

\begin{proof}
Substituting \eqref{eq:3:1:1} into \eqref{eq:2:6} and evaluating
the loop transfer function for $s = j \omega$ results in
\begin{equation}\label{eq:3:1:4}
H(j \omega) = - \frac{1}{a \, \omega^2 - c + j \omega \, (\omega^2
- b)}.
\end{equation}
Evaluating the real part of \eqref{eq:3:1:4} and substituting it
together with $k=cd$ into \eqref{eq:2:7} results in
$$
1 + \frac{c^2 d - a c d \, \omega^2}{(b \omega - \omega^3)^2 + (c
- a \, \omega^2)^2 } > 0.
$$
By comparing it with \eqref{eq:3:1:3}, the proof is completed.
\end{proof}

\vspace{1mm}

\begin{rem}
\label{rem:1} For both boundary frequencies $\omega = 0$ and
$\omega = \infty$ the condition \eqref{eq:3:1:3} is always
fulfilled. Indeed, for $\omega = \infty$ one obtains $-1 / \infty
> -1$. Also, substituting $\omega = 0$ into \eqref{eq:3:1:3} results
in $d > -1$, that is always true for the control \eqref{eq:3:1}.
\end{rem}

\begin{figure}[!h]
\centering
\includegraphics[width=0.7\columnwidth]{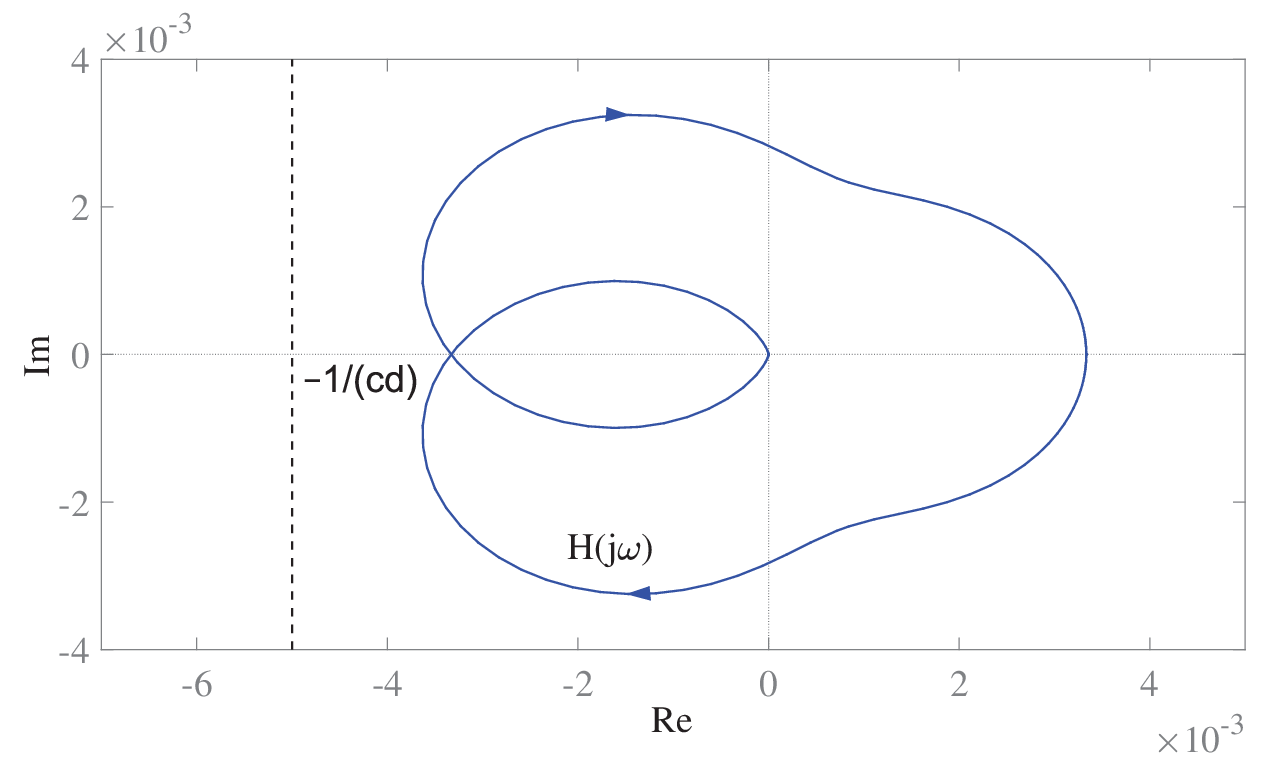}
\caption{Visualization of the circle condition for $H(j\omega)$
and sector $[0, cd]$.} \label{fig:3:0}
\end{figure}
A graphical visualization of the circle stability criterion
applied to $H(j \omega)$, with an exemplary stable characteristic
polynomial \eqref{eq:2:3} which is satisfying \eqref{eq:3:1:3} and $ k = c
d$, is depicted in Figure \ref{fig:3:0} for the sake of
a qualitative explainability.

\subsection{Disturbance rejection}
\label{sec:3:sub:2}

\begin{figure}[!h]
\centering
\includegraphics[width=0.45\columnwidth]{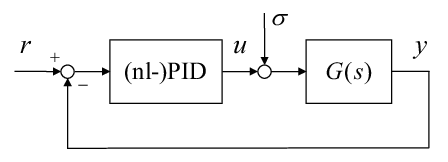}
\caption{Closed loop system with (nl-)PID controllers.}
\label{fig:3:1}
\end{figure}

Consider a standard feedback control loop as shown in Figure
\ref{fig:3:1}. Note that if a reference $r(t)\neq 0$ is also applied,
then both controllers \eqref{eq:2:2} and \eqref{eq:3:1} will
incorporate the control error $\epsilon=r-y$ instead of $y$ and,
respectively, its derivative and integral. Recall that this is a rather standard step of a coordinates transformation, since a stabilization problem (respectively disturbance rejection problem) with $r=0$ was considered before. Consequently, the control terms will
change to their additive inverse, i.e. $-a \dot{y}$ will change to $a
\dot{\epsilon}$ and so on. In the following, for the sake of
simplicity and due to the focus on disturbance rejection, we
further assume $r=0$. At the same time, an $\sigma \neq 0$ value is entering the
control loop. Further, the control \eqref{eq:3:1} will be denoted by
nl-PID while its linear counterpart \eqref{eq:2:2} by PID,
respectively, cf. Figure \ref{fig:3:1}.

\begin{figure}[!h]
\centering
\includegraphics[width=0.7\columnwidth]{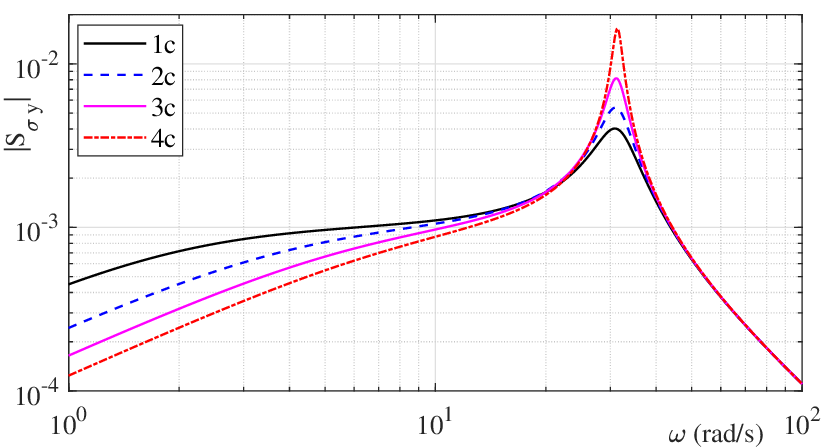}
\caption{Magnitude of frequency response function $S_{\sigma y}$
depending on $\Omega$.} \label{fig:3:2}
\end{figure}

While the transfer function of the system plant is given by
\begin{equation}\label{eq:3:2:1}
G(s) = \frac{1}{s^2},
\end{equation}
one can temporary assume \eqref{eq:3:2} as a 'frozen' (i.e.
constant) parameter and write the control transfer function as
\begin{equation}\label{eq:3:2:2}
\{\hbox{PID},\, \hbox{nl-PID} \} \equiv  C(s) = \frac{a s^2 + b s
+ \Omega}{s}.
\end{equation}
The corresponding perturbation-to-output transfer function is
\begin{equation}\label{eq:3:2:3}
S_{\sigma y}(s) = \frac{y(s)}{\sigma(s)} =\frac{G(s)}{1+C(s)G(s)}.
\end{equation}
Substituting \eqref{eq:3:2:1}, \eqref{eq:3:2:2} and evaluating
\eqref{eq:3:2:3}, one can recognize that for any $\sigma =
\mathrm{const}$, the resulting $y(s) \rightarrow 0$ as $s
\rightarrow 0$, i.e. for the steady-state. This is irrespectively of
the $\Omega$-value for which the characteristic polynomial
\eqref{eq:2:3} with $c = \Omega$ is maintained to remain stable.
Evaluating the magnitude of the frequency response function
$S_{\sigma y}(j\omega)$ one obtains
\begin{equation}\label{eq:3:2:4}
\bigl| S_{\sigma y}(j\omega) \bigr| = \frac{\omega}{ \sqrt{\bigl(b
\omega - \omega^3 \bigr)^2 + \bigl(\Omega - a \omega^2 \bigr)^2} }
\, .
\end{equation}
One can show that for an increasing $\Omega$, while satisfying
$\Omega < ab$, the peak value of \eqref{eq:3:2:4} also increases,
i.e. $\Omega \nearrow \: \Rightarrow \: \max \, \bigl| S_{\sigma
y}(j\omega) \bigr| \nearrow$. A qualitative behavior of
\eqref{eq:3:2:4} is shown in Figure \ref{fig:3:2} for the
exemplary assigned $\Omega \in \{1c, 2c, 3c, 4c \}$.

At the same time, one can also show that for frequencies lower than the
peak frequency $\omega_{\mathrm{peak}}$, an increasing $\Omega$
leads to a decreasing magnitude value, i.e. $\Omega \nearrow \:
\Rightarrow \: \bigl| S_{\sigma y}(j\omega) \bigr| \searrow$ for
$\omega \ll \omega_{\mathrm{peak}}$, cf. Figure \ref{fig:3:2}.
Here we recall that higher peak values give rise to higher
overshoots during the transients, like for example in case when
all frequencies will be excited in case of a stepwise $\sigma$
disturbance. On the contrary, lower $|S_{\sigma y}|$ values reduce
the ultimately bounded output error, correspondingly, improve the
convergence performance. Thus, an appropriately dimensioned 
$\Omega$-value, cf. \eqref{eq:3:2}, can improve both the transient
and steady-state behavior of the control system. We should
consciously recall, however, that the behavior \eqref{eq:3:2:4},
cf. also Figure \ref{fig:3:2}, discloses only a steady-state
response of the feedback control system.

\begin{figure}[!h]
\centering
\includegraphics[width=0.7\columnwidth]{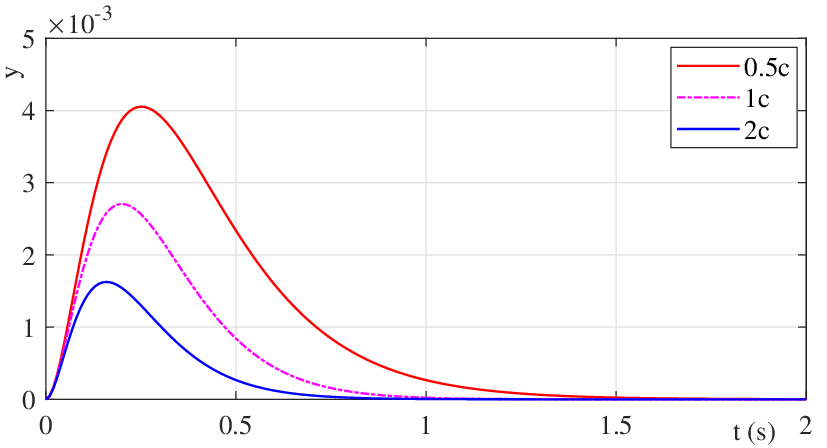}
\caption{Exemplary transient response of \eqref{eq:3:2:8} for
varying $\Omega$.} \label{fig:3:3}
\end{figure}

In order to assess the transient response of the control system to
a matched perturbation $\sigma$, and that depending on the
$\Omega$ parameter, we apply the unit step $\sigma(s)$ so that
\begin{equation}\label{eq:3:2:5}
y(s) =\frac{G(s)}{1+C(s)G(s)} \cdot \frac{1}{s}.
\end{equation}
In the following, for the sake of simplicity of developing and
analyzing the time domain solution of \eqref{eq:3:2:5}, we assume
two real eigenvalues $\lambda_1, \lambda_2 < 0$ so that the stable
transfer function \eqref{eq:3:2:3} has a double-pole pair at
$-\lambda_1$ and, respectively, one pole at $-\lambda_2$. One can
show that such poles configuration corresponds to the control
gains equal to
\begin{equation}\label{eq:3:2:6}
a = 2 \lambda_1 + \lambda_2, \quad b = \lambda_1^2 + 2 \lambda_1
\lambda_2, \quad \Omega = \lambda_1^2 \lambda_2.
\end{equation}
Solving \eqref{eq:3:2:5} in time domain and putting the common term $\Gamma = (\lambda_1 - \lambda_2)^{-2}$ outside the brackets results in
\begin{equation}\label{eq:3:2:7}
y(t) = \Gamma \Bigl(\exp(-\lambda_2 t) - \bigl(1 -
(\lambda_2-\lambda_1)\, t \bigr) \exp(-\lambda_1 t)  \Bigr).
\end{equation}
Expressing $\lambda_2$ in terms of $\Omega$, cf. \eqref{eq:3:2:6},
one obtains from \eqref{eq:3:2:7}
\begin{eqnarray}
\label{eq:3:2:8}
y(t) &=& \hat{\Gamma} \biggl(\exp\Bigl(-\frac{\Omega}{\lambda_1^2}
t\Bigr) - \Bigl( 1 - \frac{\Omega - \lambda_1^3}{\lambda_1^2} \, t
\Bigr) \exp \bigl(-\lambda_1 t\bigr) \biggr), \\
\nonumber  \hat{\Gamma} &=& \frac{\lambda_1^4}{(\lambda_1^3 -
\Omega)^2}.
\end{eqnarray}
Analyzing \eqref{eq:3:2:8}, one can show that an increasing
$\Omega$ (this for a fixed $\lambda_1$ value) leads to a decreasing
peak in the transient response of $y(t)$, cf. Figure \ref{fig:3:3}
for the exemplary taken $\Omega \in \{0.5c, 1c, 2c\}$.

\subsection*{Numerical example}
\label{sec:3:sub:3}

The following numerical example visualizes the convergence of the
control system \eqref{eq:2:1}, \eqref{eq:3:1} for various
parameter settings. The assigned (most simple) disturbance is a
constant value $\sigma=-100$. The set linear control parameters
are $a=60$, $b=1100$, $c=3000$. This is resulting in a poles
configuration $\lambda \in \{-3.28, \, -28.36 \pm j 10.47 \}$ when the nonlinear integral gain is excluded, i.e. $d=0$. Note that for an integral gain equal to $2c=6000$, that
corresponds to $d=1$ when $y \rightarrow 0$, the poles
configuration results in $\lambda^* \in \{-10, \, -20, \ -30 \}$, thus without transient oscillations in the controlled response. The initial value $y(0) = -1$ is assigned for a simultaneous set-point and disturbance rejection control task.

\begin{figure}[!h]
\centering
\includegraphics[width=0.7\columnwidth]{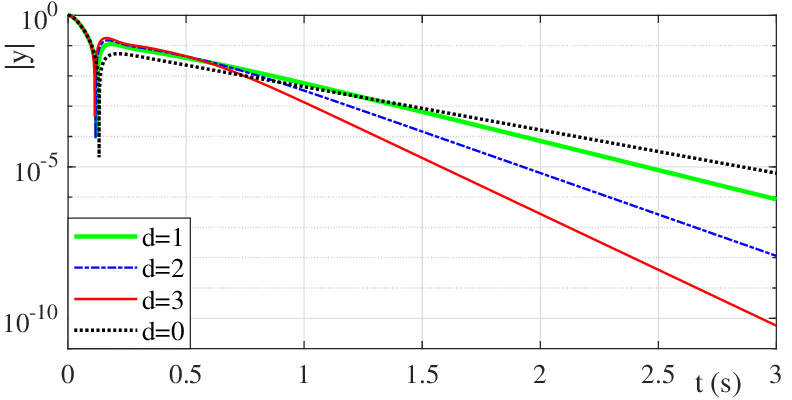}
\caption{Output convergence on logarithmic scale for varying $d$.}
\label{fig:3:2:1}
\end{figure}

For the assigned parameter $e=-10$, the output convergence on the
logarithmic scale are compared for the set parameter values $d =
\{1, 2, 3\}$ in Figure \ref{fig:3:2:1}. The linear PID case, i.e.
for $d=0$, is also included for the sake of comparison.

\begin{figure}[!h]
\centering
\includegraphics[width=0.7\columnwidth]{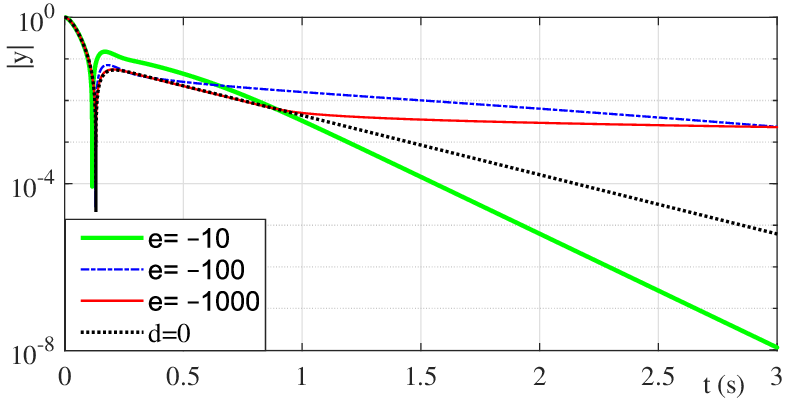}
\caption{Output convergence on logarithmic scale for varying $e$.}
\label{fig:3:2:2}
\end{figure}

Next, for the fixed parameter $d=2$, the output convergence on the
logarithmic scale is compared for the set parameter values $e =
\{-10, -100, -1000\}$ in Figure \ref{fig:3:2:2}. Also here, the
linear PID case of $d=0$ is also shown for the sake of comparison.
One can recognize that for some lower threshold $\varepsilon =
\mathrm{const}$ of the residual output error, the convergence of the
nl-PID control outperforms that one of the PID, in terms of
$|y(t)| < \varepsilon$ for $t > t_0$ with $t_0$ to be a metric of
the settling time. Simultaneously, for higher $\varepsilon$
values, for instance with $t_0 < 1 $ sec cf. Figures
\ref{fig:3:2:1}, \ref{fig:3:2:2}, it is a combination of the
$(d,e)$ parameters that allows shaping the output convergence to
be superior comparing to the linear control case. Note that in
addition to the $y(t)$-convergence, the transient over- or
undershoots and the control effort can appear as essential performance metrics for a
particular application at hand. A further explicit analysis of this, is
however out of scope of the recent communication.

\subsection{Controller design}
\label{sec:3:sub:4}

For use of the nl-PID control developed above, the controller design recommendations are summarized as follows.
\begin{enumerate}
  \item Design a standard PID feedback controller \eqref{eq:2:2} with respect to the application-specific reference following and disturbance rejecting performance criteria, see e.g. \cite{franklin2019} \cite{aastrom2021}. For the simple dynamic processes, the analytic tuning rules can be found e.g. in \cite{skogestad2003}. For an optimal disturbance rejection in the stable second-order plants by using a disturbance sensitivity function criteria, the procedure proposed in \cite{ruderman2022disturbance} can for example be used. For a detailed treatment of the PID control design, the seminal textbook e.g. \cite{aastrom2006} is recommended to be consulted.   
  \item Using the Lipschitz continuity of the perturbations, i.e. the $\Sigma$ upper bound, and the corresponding frequency characteristics of the disturbance sensitivity function \eqref{eq:3:2:4}, determine the application suitable $\Omega$-value, cf. Figure \ref{fig:3:2}.
  \item Adjust the $\Omega$-value, if necessary, in order to fulfill the transient response criteria according to \eqref{eq:3:2:7}.   
  \item For the determined $\Omega$-value, assign the $d$ parameter so that to comply with the property (i), given by \eqref{eq:3:2}, and with the stability criteria \eqref{eq:3:1:3} of Theorem \ref{theorem:2}. 
  \item Adjust the $e$ parameter so that to achieve the desired convergence behavior for a given application.
  \item Repeat the steps 2.--5. upon the numerical and experimental control evaluation if necessary.    
\end{enumerate}

\subsection{Extension by an actuator dynamics}
\label{sec:3:sub:5}

When the second-order system plant \eqref{eq:2:1} is extended by a
(parasitic) actuator dynamics
\begin{equation}\label{eq:3:3:1}
v(s) = F(s) u(s),
\end{equation}
where $F(s)$ is a stable and strictly proper transfer function
with $F(0) \rightarrow \mathrm{const}$, the system \eqref{eq:2:1}
expands to
\begin{equation}\label{eq:3:3:2}
\ddot{y}(t) = v(t) + \sigma(t).
\end{equation}
Then, the extended plant \eqref{eq:3:3:2}, \eqref{eq:3:3:1} is
controlled in the same way by \eqref{eq:3:1}, provided the
resulted Lur'e system fulfills equally the necessary assumptions
i)-iv), cf. with section \ref{sec:2:sub:2}.

For the extended plant \eqref{eq:3:3:2}, \eqref{eq:3:3:1}, the
corresponding input-to-output loop transfer function results in
\begin{equation}\label{eq:3:3:3}
H(s) = \frac{F(s)}{s^3 + a s^2 + b s +c}.
\end{equation}
Since neither the structure of a feedback system nor the
nonlinearity \eqref{eq:3:1:2} change through \eqref{eq:3:3:1}, the
same circle criterion as given in section \ref{sec:2:sub:2} is
applicable to \eqref{eq:3:3:3}. Respectively, the parametric
inequality \eqref{eq:3:1:3} of the Theorem \ref{theorem:2} should
modify, so that to take into account also the polynomial
coefficients of $F(j \omega)$. Since both the absolute stability
and disturbance rejection for \eqref{eq:3:3:3} follow the same
lines of argumentation as in sections \ref{sec:3:sub:1},
\ref{sec:3:sub:2}, and a particular actuator dynamics $F(s)$ is
rather application-specific, we omit a further detailed
development here. Also note that the controller design procedure, summarized in form of the step by step recommendations given above  
in section \ref{sec:3:sub:4}, is correspondingly extendable by the actuator dynamics $F(s)$. Here we recall that the actuator dynamics $F(s)$ will generally reduce the characteristic stability margins of the closed-loop (cf. e.g. \cite{aastrom2006,aastrom2021} for details), and thus require a PID control design to be more conservative regarding the feedback gain values.

At the same time, it is worst emphasizing that the experimental case study provided below in section \ref{sec:4} has the second-order system plant affected additionally by such
first-order actuator dynamics $v(s) = \kappa (s \mu + 1)^{-1} u(s)$, with an actuator gain factor $\kappa$ and a (not negligible) time constant $\mu > 0$.

%%%%%%%%%%%%%%%%%%%%%%%%%%%%%%%%%%%%%%%%%%%%%%%%%%%%%%%%%%%%%%%%%%%%%%%%%%%%%%%%
\section{EXPERIMENTAL EVALUATION}
\label{sec:4}

The proposed extended nonlinear PID control given by
\eqref{eq:3:1} is experimentally evaluated and additionally
benchmarked with a standard PID control \eqref{eq:2:2}. In
addition, the corresponding PD control (i.e. $c = 0$) is also
evaluated for highlighting the level of perturbations in the
system at hand and, this way, disclosing necessity of the integral
control action. The common control parameters are assigned to be the
same, while all parameter values are summarized in Table \ref{tab:1}. Note that all control gains are normalized by the overall gain $\kappa$ of the double integrator plant, cf.
\eqref{eq:3:2:1}, so that $a^* = a/\kappa$, $b^* = b/\kappa$, $c^*
= c/\kappa$, $d^* = d/\kappa$.
\begin{table}[!h]
  \renewcommand{\arraystretch}{1.4}
  \caption{Control parameters}
  \label{tab:1}
  \footnotesize
  \begin{center}
  \begin{tabular} {|p{2cm}||p{0.7cm}|p{0.7cm}|p{0.7cm}|p{0.7cm}|p{0.7cm}|}
  \hline
  control \textbackslash \hspace{0.1mm} params.   &  $a^*$ & $b^*$ & $c^*$ & $d^*$ & $e$ \\
  \hline \hline
  PD              & $10$    &  $1000$    &        &       &        \\
  \hline
  PID             & $10$    &  $1000$    & $5000$   &     &        \\
  \hline
  nl-PID          & $10$    &  $1000$    & $5000$   & $1$   &   $-100$      \\
  \hline
  \end{tabular}
  \end{center}
\end{table}

\begin{figure}[!h]
\centering
\includegraphics[width=0.3\columnwidth]{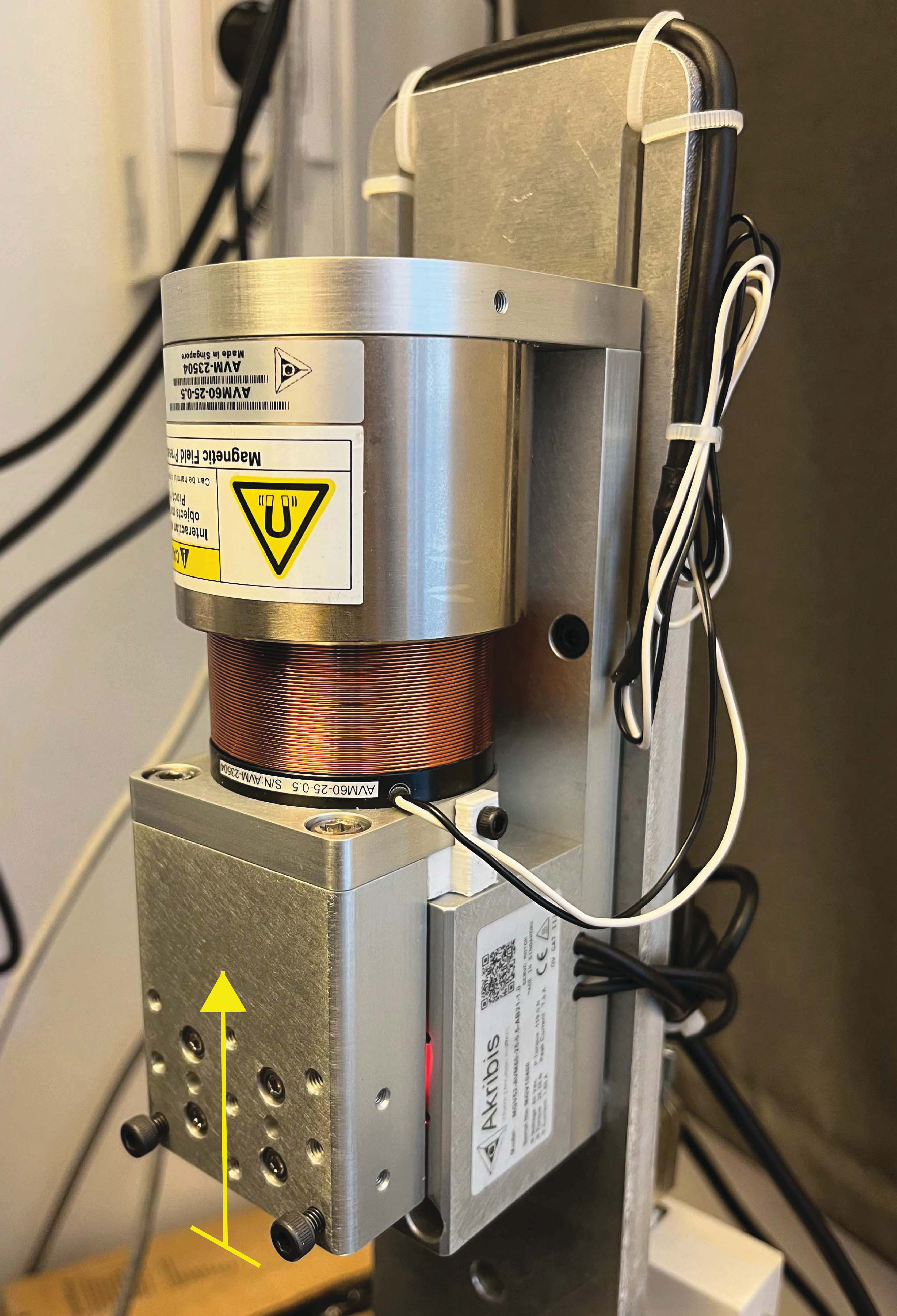}
\caption{Experimental setup of 1DOF actuator system (laboratory
view).} \label{fig:4:1}
\end{figure}

The experimental laboratory setup used for evaluation of all three
controllers is depicted in Figure \ref{fig:4:1}.
The motion system consists of a voice-coil driven actuator with
one translational degree of freedom, which has the contactless measured output position $y(t)$.
The latter contains a relatively hight level of the sensing noise.
The robust finite-time second-order sliding-mode differentiator
\cite{levant1998,moreno2012} is used for obtaining $\dot{y}(t)$.
The sampling frequency of the implemented real-time control is set
to 10 kHz. A more detailed description of the experimental system,
including the identified system parameters, can be found in
\cite{ruderman2022motion}.

Further it is worth noting that: i) an additional first-order
actuator dynamics (due to electro-magnetic behavior of the
voice-coil motor) is in place and ii) the unknown perturbation
$\sigma(t)$ is mainly a combination of the gravity, nonlinear
friction, and force ripples of the anchor-stator interaction.

\begin{figure}[!h]
\centering
\includegraphics[width=0.7\columnwidth]{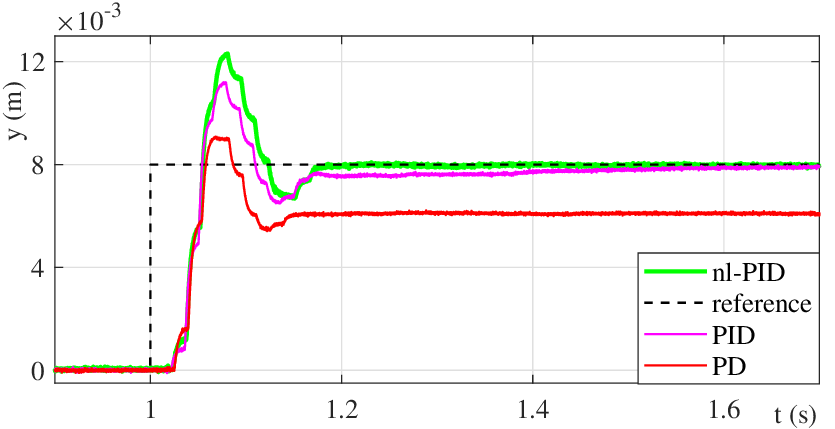}
\caption{Measured output response of the nl-PID, PID, and PD
controls.} \label{fig:4:2}
\end{figure}

\begin{figure}[!h]
\centering
\includegraphics[width=0.7\columnwidth]{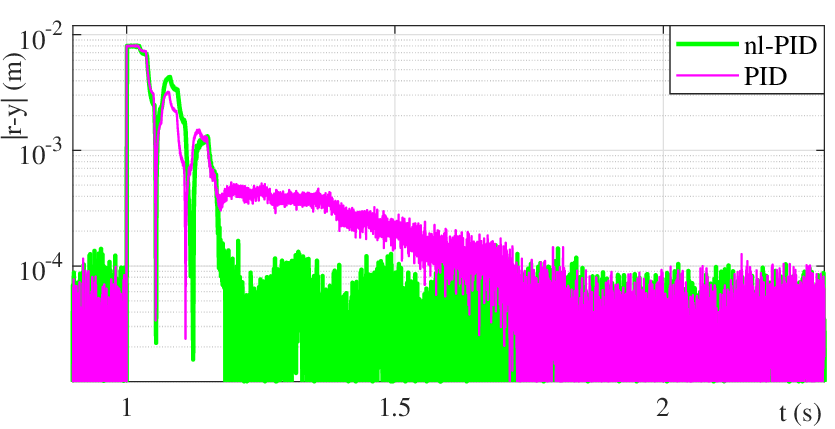}
\caption{Absolute output error (logarithmic) of the nl-PID, PID
controls.} \label{fig:4:3}
\end{figure}

The measured output response to the step reference $r = 0.008$ m,
that is applied at the time $t=1$ sec, is shown in Figure
\ref{fig:4:2} for the nl-PID, PID, and PD controllers.
The absolute output error $|\epsilon| \equiv |r-y|$ of the nl-PID
and PID controllers are additionally visualized on the logarithmic
scale in Figure \ref{fig:4:3}. One can recognize that for
$|\epsilon| < 0.0001$ m, both controllers have to deal with the
sensing noise so that the residual control errors are similar. At
the same time, the convergence performance is well distinguishable
for the time interval between 1.18 and 1.7 sec. The transient
behavior (including over- and undershoots) at the time between 1
and 1.18 sec is well comparable.

%%%%%%%%%%%%%%%%%%%%%%%%%%%%%%%%%%%%%%%%%%%%%%%%%%%%%%%%%%%%%%%%%%%%%%%%%%%%%%%%
\section{CONCLUSIONS}
\label{sec:5}

In this paper, a novel nonlinear extension of the PID feedback
control (denoted by nl-PID) is introduced. The nonlinear extension
is proposed for integral control action that allows improving
convergence performance in presence of the matched unknown
perturbations in the second-order systems. The nonlinear integral
control part does not change the system structure and is
continuous in time and Lipschitz in the system output variable, so
that the absolute stability by means of the circle criterion
turned out applicable and sufficient for the analysis. The nl-PID
control has five design parameters in total, while keeping the
same meaning and tuning rules for three of them -- proportional, integral, and
derivative feedback gains. The global asymptotic stability of the
nl-PID control was shown for the case of constant perturbations.
For the case of Lipschitz perturbations, an ultimately bounded
output error can be guaranteed by consideration in frequency
domain. Beyond the developed analysis and numerical illustrative examples, an
elaborated experimental case study is demonstrated for benchmarking the
standard PD, PID, and proposed nl-PID controllers.

%\nocite{*}% Show all bib entries - both cited and uncited; comment this line to view only cited bib entries;

\bibliographystyle{wileyNJD-AMS}
\bibliography{references}

%\section*{Author Biography}
%\begin{biography}{\includegraphics[width=60pt,height=70pt,draft]{author}}{\textbf{Michael Ruderman}
%earned his Dr.-Ing. degree in electrical engineering from TU
%University Dortmund, Germany, in 2012. He is a full professor at
%the University of Agder, Grimstad, Norway, teaching control theory
%in B.Sc., M.Sc., and Ph.D. degree programs. He serves in different
%editorial boards and technical committees of IEEE and IFAC
%societies and is chairing IEEE/IES TC on Motion Control in the
%terms 2018-2019 and 2020-2021. He is a Senior Member of IEEE and
%was the general chair of the 16th IEEE International Workshop on
%Advanced Motion Control, in 2020. His current research interests
%are in motion control, nonlinear dynamics, and hybrid control
%systems.}
%\end{biography}

\end{document}